\newcommand{\al}{\alpha}               \newcommand{\be}{\beta}
\newcommand{\ga}{\gamma}               
               \newcommand{\De}{\Delta}
\newcommand{\lb}{\lambda}              
\newcommand{\sig}{\sigma}              
               
\newcommand{\veps}{\varepsilon}        

\newcommand{\cal}{\mathcal}

\newcommand{\calv}{{\cal V}}

\newcommand{\fix}{{\rm fix}}

        \newcommand{\limpl}{\Longrightarrow}
    \newcommand{\lequi}{\Longleftrightarrow}
\newcommand{\oo}{\infty}

             \newcommand{\sk}{\smallskip}

                \def\R+oo{R_+\cup\{\oo\}}

\def\dtends   {\stackrel {\it d}{\longrightarrow}}

\def\Detends {\stackrel {\it \De}{\longrightarrow}}

\def\(V)tends {\stackrel {(\calv)}{\longrightarrow}}

\newcommand{\barr}{\begin{array}}          \newcommand{\earr}{\end{array}}
\newcommand{\beq}{\begin{equation}}        \newcommand{\eeq}{\end{equation}}
\newcommand{\bexample}{\begin{example}}    \newcommand{\eexample}{\end{example}}
\newcommand{\bit}{\begin{itemize}}         \newcommand{\eit}{\end{itemize}}
\newcommand{\blemma}{\begin{lemma}}        \newcommand{\elemma}{\end{lemma}}
\newcommand{\bprop}{\begin{proposition}}   \newcommand{\eprop}{\end{proposition}}
\newcommand{\bproof}{\begin{proof}}        \newcommand{\eproof}{\end{proof}}
\newcommand{\bremark}{\begin{remark}}      \newcommand{\eremark}{\end{remark}}
\newcommand{\btab}{\begin{tabular}}        \newcommand{\etab}{\end{tabular}}
\newcommand{\btheorem}{\begin{theorem}}    \newcommand{\etheorem}{\end{theorem}}

\documentclass[reqno]{amsart}

\newtheorem{theorem}{\bf Theorem}

\newtheorem{example}{\bf Example}
\newtheorem{lemma}{\bf Lemma}
\newtheorem{proposition}{\bf Proposition}
\newtheorem{remark}{\bf Remark}

\begin{document}

\title
[Product Fixed Points in Ordered Metric Spaces]
{PRODUCT FIXED POINTS IN \\
ORDERED METRIC SPACES}

\author{Mihai Turinici}
\address{
"A. Myller" Mathematical Seminar;
"A. I. Cuza" University;
700506 Ia\c{s}i, Romania
}
\email{mturi@uaic.ro}




\begin{abstract}
All product fixed point results 
in ordered metric spaces 
based on linear contractive conditions
are but a vectorial form 
of the fixed point statement due to
Nieto and Rodriguez-Lopez 
[Order, 22 (2005), 223-239],
under the lines in 
Matkowski 
[Bull. Acad. Pol. Sci. (Ser. Sci. Math. Astronom. Phys.), 21 (1973), 323-324].
\end{abstract}

\maketitle

\section{Introduction}
\setcounter{equation}{0}

Let $(X,d;\le)$ be a partially ordered metric space; and
$T:X\to X$ be a selfmap of $X$, with
\bit
\item[(a01)]
$X(T,\le):=\{x\in X; x\le Tx\}$ is nonempty
\item[(a02)]
$T$ is increasing ($x\le y$ implies $Tx\le Ty$).
\eit
We say that $x\in X(T,\le)$ is a {\it Picard point} (modulo $(d,\le;T)$) if
{\bf i)} $(T^nx; n\ge 0)$ is $d$-convergent, 
{\bf ii)} $z:=\lim_nT^nx$ belongs to $\fix(T)$
(in the sense: $z=Tz$).
If this happens for each $x\in X(T,\le)$,
then $T$ is referred to as a {\it Picard operator} (modulo $(d,\le)$); 
moreover, if these conditions hold for each $x\in X$, and
{\bf iii)} $\fix(T)$ is a singleton, 
then $T$ is called a {\it strong Picard operator} (modulo $(d,\le)$); 
cf. Rus \cite[Ch 2, Sect 2.2]{rus-2001}.
Sufficient conditions for such properties
are obtainable under metrical contractive requirements. 
Namely, call $T$, {\it $(d,\le;\al)$-contractive} (where $\al> 0$), if
\bit
\item[(a03)]
$d(Tx,Ty)\le \al d(x,y)$,\ $\forall x,y\in X$, $x\le y$.
\eit
Let $(x_n; n\ge 0)$ be a sequence in $X$;
call it $(\le)$-{\it ascending (descending)}, 
if $x_n\le x_m$ ($x_n\ge x_m$), provided $n\le m$.
Further, let us say that $u\in X$ is an 
{\it upper (lower) bound} of this sequence, 
when $x_n\le u$ ($x_n\ge u$), $\forall n$;
if such elements exist, we say that $(x_n; n\ge 0)$
is {\it bounded above (below)}.
Finally, call $(\le)$, {\it $d$-self-closed}
when the $d$-limit of each ascending sequence is an upper bound of it.

\btheorem \label{t1}
Assume that $d$ is complete 
and $T$ is $(d,\le;\al)$-contractive, for some $\al\in ]0,1[$.
If, in addition,
\bit
\item[(a04)]
either $T$ is $d$-continuous or $(\le)$ is $d$-self-closed
\eit
then, $T$ is a Picard operator (modulo $(d,\le)$).
Moreover, if (in addition to (a04))
\bit
\item[(a05)]
for each $x,y\in X$, $\{x,y\}$ has a lower and upper bound
\eit
then, $T$ is a strong Picard one (modulo $(d,\le)$).
\etheorem

Note that the former conclusion was obtained in 2005 by
Nieto and Rodriguez-Lopez \cite{nieto-rodriguez-lopez-2005};
and the latter one is just the 2004 result in 
Ran and Reurings \cite{ran-reurings-2004}.
For appropriate extensions of these, we refer to
Section 3 below.

According to certain authors
(cf. 
\cite{o-regan-petrusel-2008}
and the references therein)
these two results are the first extension of the 
Banach's contraction mapping principle 
to the realm of (partially) ordered metric spaces.
However, the assertion is not entirely true:
some early statements of this type 
have been obtained in 1986 by
Turinici \cite[Sect 2]{turinici-1986},
in the context of ordered metrizable uniform spaces.
\sk

Now, these fixed point results found some useful applications to 
matrix and differential/integral equations theory;
see the quoted papers for details.
As a consequence, Theorem \ref{t1} 
was the subject of many extensions.
Among these, we mention the
coupled and tripled fixed point results 
in product ordered metric spaces,
constructed under the lines in 
Bkaskar and Lakshmikantham \cite{bhaskar-lakshmikantham-2006}.
It is our aim in the following to show that, 
for all such results based on "linear" contractive conditions, 
a reduction to Theorem \ref{t1} is possible;
we refer to Section 4 and Section 5 for details.
The basic tool is the concept of {\it normal matrix} due to
Matkowski \cite{matkowski-1973}
(cf. Section 2).
Further aspects will be delineated elsewhere.

\section{Normal matrices}
\setcounter{equation}{0}

Let $R^n$ denote the usual vector $n$-dimensional space, 
$R_+^n$ the standard positive cone in $R^n$, and $\le$, the induced ordering. 
Also, let $(R_+^0)^n$ denote the interior of $R^n$ 
and $<$ the strict (irreflexive transitive) 
ordering induced by it, in the sense
\bit
\item[]
$(\xi_1,...\xi_n)< (\eta_1,...,\eta_n)$ 
provided $\xi_i< \eta_i$,  $i\in \{1,..., n\}$.
\eit
We shall indicate by $L(R^n)$ the (linear) space of all (real) $n\times n$  
matrices $A=(a_{ij})$ and by $L_+(R^n)$ the positive cone of $L(R^n)$ 
consisting of all matrices $A=(a_{ij})$ with $a_{ij}\ge 0$, $i,j\in \{1,...,n\}$. 
For each $A\in L_+(R^n)$, let us put
\bit
\item[]
$\nu(A)=\inf\{\lb\ge 0; Az\le \lb z,\ \mbox{for some}\ z> 0\}$;
\eit
and call $A$, {\it normal}, if $\nu(A)< 1$; 
or, equivalently, when the system of inequalities
\beq \label{201}
a_{i1}\zeta_1+...+a_{in}\zeta_n< \zeta_i,\ \   i\in \{1,...,n\}
\eeq
has a solution $z=(\zeta_1,...,\zeta_n)> 0$.
Concerning the problem of characterizing this class of matrices,
the following result obtained by 
Matkowski \cite{matkowski-1973}
must be taken into consideration. 
Denote (for $1\le i,j\le n$)
\bit
\item[(b01)]
$a_{ij}^{(1)}=1-a_{ij}$ if $i=j$; $a_{ij}^{(1)}=a_{ij}$, if  $i\ne j$; 
\eit
and, inductively (for $1\le k\le n-1$, $k+1\le i,j\le n$)
\bit
\item[(b02)]
$a_{ij}^{(k+1)}=a_{kk}^{(k)}a_{ij}^{(k)}-a_{ik}^{(k)}a_{kj}^{(k)}$, if $i=j$;\
$a_{ij}^{(k+1)}=a_{kk}^{(k)}a_{ij}^{(k)}+a_{ik}^{(k)}a_{kj}^{(k)}$, if $i\ne j$.
\eit

\bprop \label{p1}
The matrix $A\in L_+(R^n)$ is normal, if and only if
\bit
\item[(b03)]
$a_{ii}^{(i)}> 0$,\ \  $i\in \{1,...,n\}$.
\eit
\eprop

\bproof
{\it Necessity}. Assume that (\ref{201}) has a solution	
$z=(\zeta_1,...,\zeta_n)> 0$; that is
\beq \label{202}
\barr{rrrrrr}
a_{11}^{(1)}\zeta_1&-a_{12}^{(1)}\zeta_2&-a_{13}^{(1)}\zeta_3&-...&-a_{1n}^{(1)}\zeta_n&>0 \\
-a_{21}^{(1)}\zeta_1&+a_{22}^{(1)}\zeta_2&-a_{23}^{(1)}\zeta_3&-...&-a_{2n}^{(1)}\zeta_n&>0 \\
-a_{31}^{(1)}\zeta_1&-a_{32}^{(1)}\zeta_2&+a_{33}^{(1)}\zeta_3&-...&-a_{3n}^{(1)}\zeta_n&> 0 \\
...&... &... &...&...&... \\
-a_{n1}^{(1)}\zeta_1&-a_{n2}^{(1)}\zeta_2&-a_{n3}^{(1)}\zeta_3&-...&+a_{nn}^{(1)}\zeta_n&>0. \\
\earr
\eeq
In view of
\beq \label{203}
a_{ij}^{(1)}\ge 0,\ \ i,j\in \{1,...,n\},\ \  i\ne j,
\eeq
we must have (by these conditions)
$a_{11}^{(1)},...,a_{nn}^{(1)} > 0$;
hence, in particular, (b03) is fulfilled for $i=1$. 
Further, let us multiply the first inequality of (\ref{202}) 
by the factor $a_{i1}^{(1)}/a_{11}^{(1)}\ge 0$ 
and add it to the $i$-th relation of the same system for $i\in \{2,...,n\}$;
one gets [if we take into account (b02) (for $k=1$)
plus $a_{11}^{(1)}> 0$]
\beq \label{204}
\barr{rrrrrr}
a_{11}^{(1)}\zeta_1&-a_{12}^{(1)}\zeta_2&-a_{13}^{(1)}\zeta_3&-...&-a_{1n}^{(1)}\zeta_n&>0 \\
&a_{22}^{(2)}\zeta_2&-a_{23}^{(2)}\zeta_3&-...&-a_{2n}^{(2)}\zeta_n&>0 \\
&-a_{32}^{(2)}\zeta_2&+a_{33}^{(2)}\zeta_3&-...&-a_{3n}^{(2)}\zeta_n&> 0 \\
&... &... &...&...&... \\
&-a_{n2}^{(2)}\zeta_2&-a_{n3}^{(2)}\zeta_3&-...&+a_{nn}^{(2)}\zeta_n&>0. \\
\earr
\eeq
Since (see above)
\beq \label{205}
a_{ij}^{(2)}\ge 0,\ \  i,j\in \{2,...,n\},\ \  i\ne j,
\eeq
we must have (by these conditions)
$a_{22}^{(2)},...,a_{nn}^{(2)} > 0$;
wherefrom, (b03) is fulfilled for $i\in\{1,2\}$. 
Now, if we multiply the second inequality of (\ref{204}) 
by the factor $a_{i2}^{(2)}/a_{22}^{(2)}\ge 0$ and add it to the $i$-th relation 
of the same system for $i\in \{3,...,n\}$, one obtains that (b03) will be fulfilled 
with $i\in \{1,2,3\}$; and so on. Continuing in this way, it is clear that, after $n$ steps, (b03)
will be entirely satisfied.

{\it Sufficiency}. Assume that (b03) holds; we must find a solution
 $z=(\zeta_1,...,\zeta_n)$ for (\ref{201}) with 
$\zeta_i> 0$, $i\in \{1,...,n\}$. To do this, let us start from the system
\beq \label{206}
\barr{rrrrrr}
a_{11}^{(1)}\zeta_1&-a_{12}^{(1)}\zeta_2&-a_{13}^{(1)}\zeta_3&-...&-a_{1n}^{(1)}\zeta_n&=\sig_1 \\
-a_{21}^{(1)}\zeta_1&+a_{22}^{(1)}\zeta_2&-a_{23}^{(1)}\zeta_3&-...&-a_{2n}^{(1)}\zeta_n&=\sig_2 \\
-a_{31}^{(1)}\zeta_1&-a_{32}^{(1)}\zeta_2&+a_{33}^{(1)}\zeta_3&-...&-a_{3n}^{(1)}\zeta_n&=\sig_3 \\
...&... &... &...&...&... \\
-a_{n1}^{(1)}\zeta_1&-a_{n2}^{(1)}\zeta_2&-a_{n3}^{(1)}\zeta_3&-...&+a_{nn}^{(1)}\zeta_n&=\sig_n \\
\earr
\eeq 
where $y=(\sig_1,...,\sig_n)> 0$ is arbitrary fixed. Denote
\bit
\item[(b04)]
$\sig_i^{(1)}=\sig_i$\ (hence $\sig_i^{(1)}> 0$),\ \ $i\in \{1,...,n\}$;
\eit
and, inductively (for $1\le k\le n-1$)
\bit
\item[(b05)]
$\sig_i^{(k+1)}=a_{kk}^{(k)}\sig_i^{(k)}+a_{ik}^{(k)}\sig_k^{(k)}$,\ \ $k+1\le i\le n$.
\eit 
Let us apply to (\ref{206}) the same transformations as in (\ref{202});
one gets 
\beq \label{207}
\barr{rrrrrr}
a_{11}^{(1)}\zeta_1&-a_{12}^{(1)}\zeta_2&-a_{13}^{(1)}\zeta_3&-...&-a_{1n}^{(1)}\zeta_n&=\sig_1^{(1)} \\
&a_{22}^{(2)}\zeta_2&-a_{23}^{(2)}\zeta_3&-...&-a_{2n}^{(1)}\zeta_n&=\sig_2^{(2)} \\
&-a_{32}^{(2)}\zeta_2&+a_{33}^{(2)}\zeta_3&-...&-a_{3n}^{(2)}\zeta_n&=\sig_3^{(2)} \\
&... &... &...&...&... \\
&-a_{n2}^{(2)}\zeta_2&-a_{n3}^{(2)}\zeta_3&-...&+a_{nn}^{(2)}\zeta_n&=\sig_n^{(2)}; \\
\earr
\eeq  
where, in addition (taking (\ref{203}) into account) 
\beq \label{208}
\sig_i^{(2)}> 0,\ \  i\in \{2,...,n\}.
\eeq
If we apply to this new system the same transformations as in (\ref{204})
and, further, iterate these upon the obtained system, etc.,
we arrive at 
\beq \label{209}
\barr{rrrrrr}
a_{11}^{(1)}\zeta_1&-a_{12}^{(1)}\zeta_2&-a_{13}^{(1)}\zeta_3&-...&-a_{1n}^{(1)}\zeta_n&=\sig_1^{(1)} \\
&a_{22}^{(2)}\zeta_2&-a_{23}^{(2)}\zeta_3&-...&-a_{2n}^{(1)}\zeta_n&=\sig_2^{(2)} \\
&&a_{33}^{(3)}\zeta_3&-...&-a_{3n}^{(3)}\zeta_n&=\sig_3^{(3)} \\
&&&...&...&... \\
&&&&a_{nn}^{(n)}\zeta_n&=\sig_n^{(n)} \\
\earr
\eeq  
(the diagonal form of (\ref{206})).
From (b03) [plus the positivity properties of type (\ref{208})], 
the unique solution $z=(\zeta_1,...,\zeta_n)$ of (\ref{209})
satisfies  $\zeta_i> 0$, $i\in \{1,...,n\}$;
this, and the equivalence between (\ref{206}) and (\ref{209}), 
ends the argument.
\eproof

A useful variant of Matkowski's condition (b03) 
may now be depicted as follows. Letting $I$ denote the unitary matrix in $L(R^n)$, 
indicate by $\De_1$,...,$\De_n$ the successive "diagonal" minors of $I-A$; that is
$$
\De_1=1-a_{11},\
\De_2=\det\left(\barr{rr}1-a_{11}&-a_{12}\\-a_{21}&1-a_{22}\earr\right),...,
\De_n=\det(I-A).
$$
By the transformations we used in passing from (\ref{206} to (\ref{207}) and from this to 
the next one, etc., one gets
$\De_i=a_{ii}^{(i)},\ \ 1\le i\le n$;
so that, (b03) writes
\bit
\item[(b06)]
$\De_i> 0$,\ \ $i\in \{1,...,n\}$.
\eit
After 
Perov's terminology \cite{perov-1964}, 
a matrix $A\in L_+(R^n)$ satisfying (b06) will be termed 
{\it admissible} (or, equivalently: $a$-matrix). 
We therefore proved

\bprop \label{p2}
Over the subclass $L_+(R^n)$, we have normal $\lequi$ admissible. 
\eprop

For practical and theoretical reasons, further characterizations 
of this class of matrices are necessary.
To this end, let $||.||$  denote one of the usual 
$p$-norms in $R^n$ (where $1\le p\le \oo$),
introduced as: for $x=(\xi_1,...,\xi_n)\in R^n$,
\bit
\item[]
$||x||_p=(|\xi_1|^p+...+|\xi_n|^p)^{1/p}$ ($1\le p< \oo$);  
$||x||_\oo=\max\{|\xi_1|,...,|\xi_n|\}$.
\eit
Note that, all these have the important property
\beq \label{210}
\mbox{
$x,y\in R_+^n$, $x\le y$ $\limpl$ $||x||\le ||y||$\  \
($||.||$ is {\it monotone}).
}
\eeq
Let also $||.||^*$ stand for the compatible matrix norm in $L(R^n)$:
\bit
\item[(b07)]
$||A||^*=\sup\{||Ax||; ||x||\le 1\}$,\ \ $A\in L(R^n)$.
\eit
We stress that, for the arbitrary fixed $p\in [1,\oo]$,
the compatible to $||.||_p$ matrix norm $||.||_p^*$ in $L(R^n)$ 
is not identical with the  $p$-norm of $L(R^n)$, 
obtained by a formal identification of it with  $R^{(n^2)}$.
For example, one has
\beq \label{211}
||A||_1^*=\max\{||Ae_1||_1,...,||Ae_n||_1\},\ 
||A||_\oo^*=\max\{||e_1^\top A||_1,...,||e_n^\top A||_1\};
\eeq
where $E=(e_1,...,e_n)$ is the canonical basis of $R^n$;
hence the claim.

In the following, it will be convenient to take 
$(||.||;||.||^*)$ as $(||.||_1;||.||_1^*)$.

\bremark \label{r1}
\rm
Let $||.||_V$ be an arbitrary norm of $R^n$;
and $||.||_V^*$, its associated matrix norm in $L(R^n)$.
By a well known result
(see, for instance, 
Precupanu \cite[Ch 2, Sect 2.5]{precupanu-1992})
$||.||_V$ and $||.||$ (as defined above) are equivalent:
$$
\mbox{
$\be ||x||\le ||x||_V\le \ga ||x||$, $\forall x\in R^n$,
where $0< \be< \ga$;
}
$$
and so are the compatible matrix norms $||.||_V^*$ and $||.||^*$.
A direct verification of this last fact is to be obtained by means of 
(b07); we do not give details.
\eremark

Having these precise, call $A\in L(R^n)$, {\it asymptotic} provided 
\bit
\item[]
$A^p \to 0$ (in the matrix norm $||.||^*$) as $p\to \oo$;
\eit
or, equivalently (see (\ref{211})) if it fulfills one of the properties 
\bit
\item[]
$A^px \to 0$ as $p\to \oo$, $\forall x\in R_+^n$;\ \
$A^px \to 0$ as $p\to \oo$, $\forall x\in R^n$.
\eit
The following simple result will be in effect for us.

\blemma \label{le1}
For the matrix $A\in L_+(R^n)$ we have
\beq \label{212}
\mbox{
[$A$ is asymptotic]  $\lequi$ 
the series  $\sum_{p\ge 0} A^p$ converges in $(L(R^n),||.||^*)$. 
}
\eeq
In such a case, the sum of this series is $(I-A)^{-1}$; 
hence, $I-A$ is invertible in $L(R^n)$ and its inverse belongs to $L_+(R^n)$.
\elemma

\bproof
Let the matrix $A$ be asymptotic. If $x\in R^n$ satisfies $(I-A)x=0$ then,
(by repeatedly applying $A$ to the equivalent equality) 
$x=A^px$,\ \  for all $p\in N$;
wherefrom $x=0$ (if one takes the limit as $p\to \oo)$; proving that  
$(I-A)^{-1}$ exists as an element of $L(R^n)$. Moreover, in view of
$I-A^p=(I-A)(I+A+...+A^{p-1})$,  $p> 1$,
one gets (again by a limit process)
$I=(I -A)(I+A+A^2+...)$; which ends the argument.
\eproof

As before, we may ask of which relationships exist between 
this class of matrices and the preceding ones.
To do this, the following renorming  statement  
involving  normal  matrices will be useful.

\blemma \label{le2}
Let $A\in L_+(R^n)$ be a normal matrix. 
Then, an equivalent monotonic norm $||.||_A$ in $R^n$ 
and a number $\al$ in $]0,1[$ exist with the property
\beq \label{213}
\mbox{
$||Ax||_A\le \al ||x||_A$,\ \ for all $x \in R_+^n$. 
}
\eeq
\elemma

\bproof
By hypothesis, we have promised a vector $z=(\zeta_1,...,\zeta_n)> 0$ 
and a number $\al$ in  $]\nu(A),1[$ with $Az\le \al z$. 
Define a norm  $||.||_A$ in $R^n$ as
\bit
\item[(b08)]
$||x||_A=\max\{|\xi_i|/\zeta_i; 1\le i\le n\}$, $x=(\xi_1,...,\xi_n)\in R^n$.
\eit
By the obvious relation
\beq \label{214}
\mbox{
$x\le (||x||_A)z$,\ \  for all  $x\in R_+^n$,
}
\eeq
one gets (if we take into account the choice of $z$)
$$
Ax\le (||x||_A)Az \le (\al ||x||_A) z,\ \ x\in R_+^n;
$$
wherefrom, (\ref{213}) results. 
Since the monotonic property is evident, we omit the details. 
It remains only to prove that $||.||_A$ is equivalent with  
the initial norm $||.||$ in $R^n$. 
But this follows easily by the relation (deduced from (b08) and (\ref{214}))
\beq \label{215}
\be ||x||_A \le ||x||\le \ga ||x||_A,\ \ x\in R^n;
\eeq
where $\be=\min\{\zeta_i; 1\le i\le n\}$, $\ga=\zeta_1+...+\zeta_n$.
The proof is complete.
\eproof

We may now give an appropriate answer to the above posed problem.

\bprop \label{p3}
For each matrix of $L_+(R^n)$, we have normal $\lequi$ asymptotic. 
\eprop

\bproof
Let $A \in L_+(R^n)$ be normal. By Lemma \ref{le2}, we found an equivalent
monotonic norm $||.||_A$ on $R^n$, and an $\al\in ]0,1[$, with the property (\ref{213}). 
From this, we get
$$
\mbox{
$A^p x\to 0$ (modulo $||.||_A$) as $p \to \oo$, for all $x\in R_+^n$,
}
$$
which, according to (\ref{215}), is just the asymptotic property. 
Conversely, let $A\in L_+(R^n)$ be asymptotic. 
Fix $b=(\be_1,...,\be_n)> 0$ and put
$z=\sum_{p\ge 0} A^p b$  (hence, $z> 0$).
As $Az=\sum_{p\ge 1}A^p b$, we have $z=b+Az$; which, combined with the
choice of $b$, gives $Az< z$. 
\eproof

We cannot close these developments without giving another characterization 
of asymptotic (or normal) matrices in terms of {\it spectral radius};
this fact -- of marginal importance for the next section -- is, however,
sufficiently interesting by itself to be added here. 
Let $A\in L(R^n)$ be a matrix. Under the natural immersion 
of $R^n$ in $C^n$, let us call the number $\lb\in C$ an {\it eigenvalue} of $A$,
provided $Az=\lb z$, for some different from zero vector $z\in C^n$ 
(called in this case an {\it eigenvector} of $A$);
the class of all these numbers will be denoted 
$\sig(A)$ (the {\it spectrum} of $A$). Define
\bit
\item[]
$\rho(A)=\sup\{|\lb|; \lb\in \sig(A)\}$\ \  
(the {\it spectral radius} of $A$).
\eit

\bprop \label{p4}
The matrix $A\in L_+(R^n)$  is asymptotic if and only if $\rho(A)< 1$.
\eprop

\bproof
Suppose $A$ is asymptotic. For each eigenvalue $\lb\in \sig(A)$,
let $z=z(\lb)\in C^n$ be the corresponding eigenvector of $A$. 
We  have $Az=\lb z$; and this gives
$A^pz=\lb^p z$,\ \  for all $p\in N$.
By the choice of $A$ plus $z\ne 0$, we must have $\lb^p \to 0$ as $p\to \oo$,
which cannot happen unless  $|\lb|< 1$; hence $\rho(A)< 1$.
Conversely, assume that the matrix $A=(a_{ij})$ in $L_+(R^n)$ 
satisfies $\rho(A)< 1$; and put $A^{(\veps)}=(a_{ij}^{(\veps)})$, $\veps > 0$, where 
$$
a_{ij}^{(\veps)}=a_{ij}+\veps,\ \  1\le i,j\le n.
$$
We have $\rho(A^{(\veps)})< 1$, when $\veps > 0$ is small enough 
(one may follow a direct argument based on the obvious fact:
for each (nonempty) compact $K$ of $R$,
$$
\mbox{
$\det(\lb I-A^{(\veps)}) \to \det(\lb I-A)$ when $\veps\to 0+$,
uniformly over $\lb\in K$;
}
$$
we do not give further details).
Now, as $A^{(\veps)}$ is a matrix over $R_+^0$ 
(in the sense: $a_{ij}^{(\veps)}> 0$, $i,j\in \{1,...,n\}$),
we have, by the Perron-Frobenius theorem 
(see, e.g., Bushell \cite{bushell-1973}), 
that  for a sufficiently small  $\veps> 0$, 
$A^{(\veps)}$ has a positive eigenvalue $\mu=\mu(\veps)> 0$ 
(which, in view of $\rho(A^{(\veps)})< 1$, must satisfy $\mu< 1$), 
as well as a corresponding  eigenvector $z=z(\veps)> 0$. 
But then,
$Az\le A^{(\veps)} z= \mu z< z$; hence, $A$ is normal. 
This, along with Proposition \ref{p3}, completes the argument.
\eproof

\section{Extension of Theorem \ref{t1}}
\setcounter{equation}{0}

Let $(X,d)$ be a metric space; 
and $(\le)$ be a {it quasi-order}
(i.e.: reflexive and transitive relation) over $X$.
For each $x,y\in X$, denote:
$x<> y$ iff either $x\le y$ or $y\le x$
(i.e.: $x$ and $y$ are comparable).
This relation is reflexive and symmetric; 
but not in general transitive.
Given $x,y\in X$, any subset 
$\{z_1,...,z_k\}$ (for $k\ge 2$) in $X$ with
$z_1=x$, $z_k=y$, and [$z_i<> z_{i+1}$,  $i\in \{1,...,k-1\}$]
will be referred to as a {\it $<>$-chain} 
between $x$ and $y$; the class of all these will 
be denoted as $C(x,y;<>)$.
Let $\sim$ stand for the relation (over $X$): 
$x\sim y$ iff $C(x,y;<>)$ is nonempty. 
Clearly, $(\sim)$ is reflexive and symmetric; because so is 
$<>$. Moreover, $(\sim)$ is transitive; 
hence, it is an equivalence over $X$.
Call $d$, {\it $(\le)$-complete} when each ascending 
$d$-Cauchy sequence is $d$-convergent.
Finally, let $T:X\to X$ be a selfmap of $X$;
we say that it is $(d,\le)$-{\it continuous} when 
[$(x_n)$=ascending, $x_n\dtends x$] imply $Tx_n\dtends Tx$.

\btheorem \label{t2}
Assume (under (a01) and (a02))  that $d$ is $(\le)$-complete 
and $T$ is $(d,\le;\al)$-contractive, for some $\al\in ]0,1[$.
If, in addition,
\bit
\item[(c01)]
either $T$ is $(d,\le)$-continuous or 
$(\le)$ is $d$-self-closed,
\eit
then, $T$ is a Picard operator (modulo $(d,\le)$).
Moreover, if (in addition to (c01))
\bit
\item[(c02)]
$(\sim)=X\times X$
[$C(x,y; <>)$ is nonempty, for each $x,y\in X$],
\eit
then, $T$ is a strong Picard operator (modulo $(d,\le)$).
\etheorem

This result is a weaker form of Theorem \ref{t1};
because (a04) $\limpl$ (c01), (a05) $\limpl$ (c02).
[In fact, given $x,y\in X$,  
there exist, by (a05), some $u,v\in X$ with
$u\le x\le v$, $u\le y\le v$. This yields $x<>u$, $u<> y$; 
wherefrom, $x\sim y$].
Its proof mimics, in fact, the one of Theorem \ref{t1}.
However, for completeness reasons, we shall provide it, with
some modifications.

\bproof

{\bf I)} 
Let $x\in X(T,\le)$ be arbitrary fixed; 
and put $x_n=T^nx$, $n\in N$. 
By (a02) and (a03), 
$d(x_{n+1},x_{n+2})\le  \al d(x_n,x_{n+1}))$, for all $n$.
This yields 
$d(x_n,x_{n+1})\le \al^n d(x_0,x_1))$, $\forall n$;
so that, as the series $\sum_n \al^n$ converges, 
$(x_n; n\ge 0)$ is an ascending $d$-Cauchy sequence.
Combining with the $(\le)$-completeness of $d$,
it results that $x_n\dtends x^*$, for some $x^*\in X$.
Now, if the first half of (c01) holds, we have
$x_{n+1}=Tx_n\dtends Tx^*$; so that 
(as $d$=metric), $x^*\in \fix(T)$.
Suppose that the second half of (c01) is valid;
note that, as a consequence, $x_n\le x^*$, $\forall n$.
By the contractive condition, we derive
$d(x_{n+1},Tx^*)\le \al d(x_n,x^*)$, $\forall n$;
so that, by the obtained convergence property,
$x_{n+1}=Tx_n\dtends Tx^*$; wherefrom 
(see above) $x^*\in \fix(T)$.

{\bf II)}
Take $a,b\in X$, $a\le b$.
By the contractive condition, 
$d(T^na,T^nb)\le \al^n d(a,b)$, $\forall n$;
whence $\lim_n d(T^na,T^nb)=0$. 
From the properties of the metric, one gets 
$\lim_n d(T^na,T^nb)=0$ if $a<> b$; 
as well as (by definition)
$\lim_n d(T^na,T^nb)=0$ if $a\sim b$. 
This, along with (c02), gives the desired conclusion.
\eproof

\section{Vector linear contractions}
\setcounter{equation}{0}

Let $X$ be an abstract set; and $q\ge 1$ be a positive integer. 
In the following, the notion of $R^q$-{\it valued metric} on $X$ 
will be used to designate  any function $\De:X^2\to R_+^q$, supposed to be
{\it reflexive sufficient} [$\De(x,y)=0$ iff $x=y$]
{\it triangular} [$\De(x,z)\le \De(x,y)+\De(y,z)$, $\forall x,y,z\in X$]
and
{\it symmetric} [$\De(x,y)=\De(y,x)$, $\forall x,y\in X$].
In this case, the couple $(X,\De)$ will be termed an 
{\it $R^q$-valued metric space}.     
Fix in the following such an object; as well the usual
norm $||.||:=||.||_1$, over $R^q$. 
Note that, in such a case, the map
\bit
\item[(d01)]
($d:X^2\to R_+$):\ \ $d(x,y)=||\De(x,y)||$,\ $x,y\in X$
\eit
is a (standard) metric on $X$.
Let  also $(\preceq)$ be a {\it quasi-ordering} over $X$.

Define a $\De$-convergence property over $X$ as:
[$x_n\Detends x$ iff $\De(x_n,x)\to 0$]. 
The set of all such $x$ will be denoted $\lim_n (x_n)$; when it is
nonempty (hence, a singleton), $(x_n)$ will be termed 
{\it $\De$-convergent}.
Further, call $(x_n)$, {\it $\De$-Cauchy} provided 
[$\De(x_i,x_j)\to 0$ as $i,j\to \oo$]. 
Clearly, each $\De$-convergent sequence is $\De$-Cauchy; 
but the converse is not general valid.
Note that, in terms of the associated metric $d$, 
\beq \label{401}
\mbox{
[$\forall (x_n)$, $\forall x$]:\ $x_n\Detends x$ iff $x_n\dtends x$
}
\eeq
\beq \label{402}
\mbox{
[$\forall (x_n)$]:\ \ $(x_n)$ is $\De$-Cauchy iff $(x_n)$ is $d$-Cauchy.
}
\eeq
Call $\De$, {\it $(\preceq)$-complete} when each ascending 
$\De$-Cauchy sequence is $\De$-convergent.
Likewise, call $(\preceq)$, {\it $\De$-self-closed} when the $\De$-limit
of each ascending sequence is an upper bound of it.
By (\ref{401}) and (\ref{402}) we have the global properties
\beq \label{403}
\mbox{
[$\De$ is $(\preceq)$-complete] iff [$d$ is $(\preceq)$-complete]
}
\eeq
\beq \label{404}
\mbox{
[$(\preceq)$ is $\De$-self-closed] iff [$(\preceq)$ is $d$-self-closed].
}
\eeq
Finally, take a selfmap $T:X\to X$, according to 
\bit
\item[(d02)]
$X(T,\preceq):=\{x\in X; x\preceq Tx\}$ is nonempty
\item[(d03)]
$T$ is increasing ($x\preceq y$ implies $Tx\preceq Ty$).
\eit
We say that $x\in X(T,\preceq)$ is a {\it Picard point} 
(modulo $(\De,\preceq;T)$) if
{\bf j)} $(T^nx; n\ge 0)$ is $\De$-convergent, 
{\bf jj)} $z:=\lim_n(T^nx)$ belongs to $\fix(T)$.
If this happens for each $x\in X(T,\preceq)$,
then $T$ is referred to as a {\it Picard operator}
(modulo $(\De,\preceq)$). 
Sufficient conditions for such properties are to be 
obtained under vectorial contractive requirements.
Given  $A\in L_+(R^q)$,
let us say that $T$ is {\it $(\De,\preceq;A)$-contractive}, provided
\bit
\item[(d04)]
$\De(Tx,Ty)\le  A \De(x,y)$,\ \ $\forall x,y\in X$, $x\preceq y$.
\eit
Further, let us say that $T$ is $(\De,\preceq)$-{\it continuous} when 
[$(x_n)$=ascending and $x_n\Detends x$] imply 
$Tx_n \Detends Tx$.
As before, in terms of the associated via (d01) metric $d$, we have
(by means of (\ref{401}) and (\ref{402}) above)
\beq \label{405}
\mbox{
[$T$ is $(\De,\preceq)$-continuous] iff [$T$ is $(d,\preceq)$-continuous].
}
\eeq
The following answer to the posed question is available.

\btheorem \label{t3}
Assume (under (d02) and (d03)) that 
$\De$ is $(\preceq)$-complete 
and there exists a normal $A\in L_+(R^q)$ 
such that $T$ is $(\De,\preceq;A)$-contractive. 
In addition, suppose that
\bit
\item[(d05)]
either ($T$ is $(\De,\preceq)$-continuous) 
or 
($(\preceq)$ is $\De$-self-closed).
\eit
Then, $T$ is a Picard operator (modulo $(\De,\preceq)$).
\etheorem

\bproof 
As $A$ is normal, there exist, by Lemma \ref{le2},
an equivalent (with $||.||$) monotonic norm $||.||_A$ on $R^q$, and 
an $\al\in ]0,1[$, fulfilling (\ref{213}).
Define a new metric $e(.,.)$ over $X$, according to
\bit 
\item[(d06)]
$e(x,y)=||\De(x,y)||_A$, \ \ $x,y\in X$.
\eit
By the norm equivalence (\ref{215}), 
the properties (\ref{401})-(\ref{405})
written in terms of $d$ continue to hold in terms of $e$.
Moreover, by the monotonic property and (\ref{213}),
\beq \label{406}
e(Tx,Ty)\le \al e(x,y),\ \ \forall x,y\in X,\ x\preceq y.
\eeq
Summing up, Theorem \ref{t2} is applicable to
$(X,e;\preceq)$ and $T$; wherefrom, all is clear.
\eproof

In particular, when $(\preceq)=X^2$ 
(the {\it trivial quasi-order} on $X$)
the corresponding version of Theorem \ref{t3}
is just the statement in 
Perov \cite{perov-1964}.

\section{Product fixed points}
\setcounter{equation}{0}

Let $\{(X_i,d_i;\le_i); 1\le i\le q\}$ be a system of 
quasi-ordered metric spaces.
Denote $X=\prod\{X_i; 1\le i\le q\}$ 
(the Cartesian product of the ambient sets);
and put, for $x=(x_1,...,x_q)$ and $y=(y_1,...,y_q)$ in $X$
\bit
\item[(e01)]
$\De(x,y)=(d_1(x_1,y_1),...,d_q(x_q,y_q))$,
\item[(e02)]
$x\preceq y$ iff $x_i\le_i y_i$, $i\in \{1,...,q\}$.
\eit
Clearly, $\De$ is a $R^q$-valued metric on $X$; 
and $(\preceq)$ acts as a quasi-ordering over the same.
As a consequence of this, we may now introduce all conventions in
Section 4.
Note that, by the very definitions above, we have, 
for the sequence $(x^n=(x^n_1,...,x^n_q); n\ge 0)$ in $X$
and the point $x=(x_1,...,x_q)$ in $X$,
\beq \label{501}
\mbox{
$x^n\Detends x$ iff $d_i(x^n_i,x_i)\to 0$ as $n\to \oo$, 
for all $i\in \{1,...,q\}$
}
\eeq
\beq \label{502}
\mbox{
$(x^n; n\ge 0)$ is $\De$-Cauchy iff $(x^n_i; n\ge 0)$ is
$d_i$-Cauchy, $\forall i\in \{1,...,q\}$.
}
\eeq
This yields the useful global implications
\beq \label{503}
\mbox{
[$d_i$ is $(\le_i)$-complete, $\forall i\in \{1,...,q\}$]
$\limpl$ $\De$ is $(\preceq)$-complete
}
\eeq
\beq \label{504}
\mbox{
[$(\le_i)$ is $d_i$-self-closed, $\forall i\in \{1,...,q\}$]
$\limpl$ $(\preceq)$ is $\De$-self-closed.
}
\eeq

{\bf (I)}
We are now passing to our effective part.
Let $(T_i:X\to X_i; 1\le i\le q)$ be a system of maps;
it generates an {\it associated} selfmap (of $X$)
\bit
\item[(e03)]
$T:X\to X$:\ \  $Tx=(T_1x,...,T_qx)$, $x=(x_1,...,x_q)\in X$. 
\eit
Suppose that
\bit
\item[(e04)]
$\exists a=(a_1,...,a_q)\in X:\ a_i\le_i T_ia,\ i\in \{1,...,q\}$
\item[(e05)]
$T_i$ is increasing ($x\preceq y$ $\limpl$ $T_ix\le_iT_iy$),\ $i\in \{1,...,q\}$.
\eit
Note that, as a consequence, (d02) and (d03) hold.
For $i\in \{1,...,q\}$, call $T_i$, {\it $(\De,\preceq)$-continuous}, when:
[$(x^n=(x^n_1,...,x^n_q))$=ascending and $x^n\Detends x$] imply 
$d_i(T_ix^n,T_ix)\to 0$ as $n\to \oo$. Clearly,
\beq \label{505}
\mbox{
[$T_i$ is $(\De,\preceq)$-continuous, $i\in \{1,...,q\}$]
implies $T$ is $(\De,\preceq)$-continuous.
}
\eeq
Let $A=(a_{ij}; 1\le i,j\le q)$ be an element of $L_+(R^q)$.
For $i\in \{1,...,q\}$, denote $A_i=(a_{i1},...,a_{iq})$
(the $i$-th line of $A$).
Call  $T_i$, $(\De,\preceq;A_i)$-contractive, provided
\bit
\item[(e06)]
$d_i(T_ix,T_iy)\le  A_i \De(x,y)$,\ \ $\forall x,y\in X$, $x\preceq y$.
\eit
The following implication is evident:
\beq \label{506}
\mbox{
[$T_i$ is $(\De,\preceq;A_i)$-contractive, $i\in \{1,...,q\}$]
$\limpl$ $T$ is $(\De,\preceq;A)$-contractive.
}
\eeq
Putting these together, we have (via Theorem \ref{t3} above):

\btheorem \label{t4}
Assume (under (e04) and (e05)) that 
$d_i$ is $(\le_i)$-complete, $\forall i\in \{1,...,q\}$,
and there exists a normal matrix $A=(A_1,...,A_q)^\top \in L_+(R^q)$ 
such that $T$ is $(\De,\preceq;A_i)$-contractive, 
$\forall i\in \{1,...,q\}$. 
In addition, suppose that
\bit
\item[(e07)]
either
($T_i$ is $(\De,\preceq)$-continuous, $\forall i\in \{1,...,q\}$) \\
or ($(\le_i)$ is $d_i$-self-closed, $\forall i\in \{1,...,q\}$).
\eit
Then, the associated selfmap $T$ is a Picard one (modulo $(\De,\preceq)$).
\etheorem

In particular, when $(\le_i)=X_i\times X_i$, $i\in \{1,...,q\}$,
this result is just the one in 
Matkowski \cite{matkowski-1973}.
Some "uniform" versions of it were obtained in
Czerwik \cite{czerwik-1976}; 
see also
Balakrishna Reddy and Subrahmanyam  \cite{balakrishna-reddy-subrahmanyam-1981}. 
\sk

{\bf (II)}
By definition, any fixed point of the associated selfmap $T$ will 
be referred to as a {\it product fixed point} of 
the original system $(T_1,...,T_q)$.
To see its usefulness, it will suffice noting
that, by an appropriate choice of our data, one gets 
(concrete) coupled and tripled fixed point results in the area,
obtainable via "linear" type contractive conditions.
The most elaborated one, due to
Berinde and Borcut  \cite{berinde-borcut-2011}
will be discussed below.
\sk

Let $(X,d;\le)$ be a partially ordered metric space;
and take a map $F:X^3\to X$.
We say that $b=(b_1,b_2,b_3)\in X^3$ is 
a {\it tripled} fixed point of $F$, provided 
\bit
\item[(e08)]
$b_1=F(b_1,b_2,b_3)$, $b_2=F(b_2,b_1,b_2)$, $b_3=F(b_3,b_2,b_1)$.
\eit
Sufficient conditions for the existence of such points are centered on
\bit
\item[(e09)]
there exists at least one $a=(a_1,a_2,a_3)\in X^3$ with \\
$a_1\le F(a_1,a_2,a_3)$, $a_2\ge F(a_2,a_1,a_2)$, $a_3\le F(a_3,a_2,a_1)$
\item[(e10)]
$F$ is {\it mixed monotone}: \\
$x_1\le x_2$ $\limpl$ $F(x_1,y,z)\le F(x_2,y,z)$, \\
$y_1\le y_2$ $\limpl$ $F(x,y_1,z)\ge F(x,y_2,z)$, \\
$z_1\le z_2$ $\limpl$ $F(x,y,z_1)\le F(x,y,z_2)$.
\eit
Let $(\preceq)$ be the ordering on $X^3$ introduced as
\bit
\item[]
$(x_1,x_2,x_3)\preceq (y_1,y_2,y_3)$ iff 
$x_1\le y_1$, $x_2\ge y_2$, $x_3\le y_3$.
\eit
Call $F$, $(d,\preceq;\al_1,\al_2,\al_3)$-contractive (where $\al_1,\al_2,\al_3> 0$) when
\bit
\item[(e11)]
$d(F(x_1,x_2,x_3),F(y_1,y_2,y_3))\le 
\al_1 d(x_1,y_1)+\al_2 d(x_2,y_2)+\al_3 d(x_3,y_3)$,
for all $x=(x_1,x_2,x_3)$, $y=(y_1,y_2,y_3)$ in $X^3$ with $x\preceq y$.
\eit

\btheorem  \label{t5}
Suppose (under (e09) and (e10)) that 
$d$ is a complete metric and $F$ is 
$(d,\preceq;\al_1,\al_2,\al_3)$-contractive, for some $\al_1,\al_2,\al_3>0$
with $\al:=\al_1+\al_2+\al_3< 1$. If, in addition,
either 
[$F$ is continuous] 
or
[both $(\le)$ and $(\ge)$ are $d$-self-closed]
then $F$ has at least one tripled fixed point.
\etheorem

See the quoted paper for the original argument.
Here, we shall develop a different one, 
based on the fact that, any tripled fixed point 
for $F$ is a fixed point of the associated 
selfmap $T$ of $X^3$, introduced as:
\bit
\item[(e12)]
$Tx=(F(x_1,x_2,x_3),F(x_2,x_1,x_2),F(x_3,x_2,x_1))^\top$,
$x=(x_1,x_2,x_3)\in X^3$.
\eit
To do this, it will suffice verifying that conditions of
Theorem \ref{t4} are fulfilled with
$(X_1,d_1;\le_1)=(X,d;\le)$, 
$(X_2,d_2;\le_2)=(X,d;\ge)$,
$(X_3,d_3;\le_3)=(X,d;\le)$.

\bproof ({\bf Theorem \ref{t5}})
Define, for 
$x=(x_1,x_2,x_3)$ and $y=(y_1,y_2,y_3)$ in $X^3$,
\bit
\item[]
$\De(x,y)=(d(x_1,y_1),d(x_2,y_2),d(x_3,y_3))$.
\eit
We have to establish that the associated map $T$ 
introduced via (e12) is increasing (modulo $(\preceq)$)
and $(\De,\preceq;A)$-contractive, for a normal 
matrix $A\in L_+(R^3)$,
The former of these is directly obtainable by means of the
mixed monotone property (e10).
For the latter, note that, by (e11), $T$ is 
$(\De,\preceq;A)$-contractive, where $A\in L_+(R^3)$ is given as
$A=(A_1,A_2,A_3)^\top$, where
$$
A_1=(\al_1,\al_2,\al_3), A_2=(\al_2,\al_1+\al_3,0),
A_3=(\al_3,\al_2,\al_1).
$$
Since, on the other hand, $A\Theta=\al \Theta< \Theta$, 
where $\Theta=(1,1,1)^\top$, 
it results that $A$ is normal (cf. Section 2); and we are done.
\eproof

\bremark \label{r2}
\rm
The last part of the argument above suggests us a simplified 
proof of the original Berinde-Borcut argument. 
Namely, given $(X,d;\le)$, $F$ and $(\al_1,\al_2,\al_3)$ as in Theorem \ref{t5}, 
define a standard metric $D(.,.)$  over $X^3$ as:
for  $x=(x_1,x_2,x_3)$ and $y=(y_1,y_2,y_3)$ in $X^3$,
\bit
\item[(e13)]
$D(x,y)=\max\{d(x_1,y_1),d(x_2,y_2),d(x_3,y_3)\}$.
\eit
By the contractive condition (e11), it is clear that
\beq \label{507}
\mbox{
$D(Tx,Ty)\le \al D(x,y)$, for all $x,y\in X^3$, $x\preceq y$.
}
\eeq
This, along with the previous remarks, tells us that 
Theorem \ref{t1} applies to the ordered metric space 
$(X^3,D;\preceq)$ and  $T$; wherefrom, all is clear.
\eremark

Note, finally, that the original coupled fixed point statement in 
Bhaskar and Lakshmikantham \cite{bhaskar-lakshmikantham-2006}
corresponds to the normal matrix (where $0< \al< 1$)
$$
A=(A_1,A_2)^\top\in L_+(R^2):\ \  A_1=A_2=(\al/2)(1,1).
$$
Further aspects will be delineated elsewhere.



\begin{thebibliography}{99}


\bibitem{balakrishna-reddy-subrahmanyam-1981}
{K. Balakrishna Reddy} and {P. Subrahmanyam},
\it Extensions of Krasnoselskij's and Matkowski's fixed point theorems,
\rm Funkc. Ekv., 24 (1981),  67-83.


\bibitem{berinde-borcut-2011}
{V. Berinde} and {M. Borcut},
\it Tripled fixed point theorems for contractive type mappings 
in partially ordered metric spaces,
\rm Nonlinear Anal., 74 (2011), 4889-4897.


\bibitem{bhaskar-lakshmikantham-2006}
{T. G. Bhaskar} and {V. Lakshmikantham},
\it Fixed point theorems in partially ordered metric spaces and applications,
\rm Nonlinear Anal., 65 (2006), 1379-1393.


\bibitem{bushell-1973}
{P. J. Bushell},
\it Hilbert's metric and positive contraction mappings in a Banach space, 
\rm Arch. Rational Mech. Anal., 52 (1973), 330-338.


\bibitem{czerwik-1976}
{S. Czerwik},
\it A fixed point theorem for a system of multivalued transformations,
\rm Proc. Amer. Math. Soc., 55 (1976), 136-139.

\bibitem{matkowski-1973}
{J. Matkowski}, 
\it Some inequalities and a generalization of Banach's principle,
\rm Bull. Acad. Pol. Sci. (Ser. Sci. Math. Astronom. Phys.), 21 (1973), 323-324.


\bibitem{nieto-rodriguez-lopez-2005}
{J.J. Nieto} and  {R. Rodriguez-Lopez}, 
\it Contractive mapping theorems in partially ordered sets and applications 
to ordinary differential equations,
\rm Order, 22 (2005), 223-239.


\bibitem{o-regan-petrusel-2008}
{D. O'Regan} and  {A. Petru\c{s}el}, 
\it Fixed point theorems for generalized contractions in ordered metric spaces, 
\rm J. Math. Anal. Appl., 341 (2008), 1241-1252.


\bibitem{perov-1964}
{A. I. Perov},
\it On the Cauchy problem for systems of ordinary differential equations,
\rm (Russian), in "Approximate Methods for solving Differential Equations", 
pp. 115-134, Naukova Dumka. Kiev, 1964.


\bibitem{precupanu-1992}
{T. Precupanu},
\it Linear Topological Spaces and Fundamentals of Convex Analysis,
\rm (Romanian), Editura Academiei Rom\^{a}ne, Bucure\c{s}ti, 1992.


\bibitem{ran-reurings-2004}
{A. C. M. Ran} and  {M. C. Reurings}, 
\it A fixed point theorem in partially ordered sets and some applications to matrix equations, 
\rm Proc. Amer. Math. Soc., 132 (2004), 1435-1443.


\bibitem{rus-2001}
{I. A. Rus}, 
\it Generalized Contractions and Applications 
\rm Cluj University Press, Cluj-Napoca, 2001.


\bibitem{turinici-1986}
{M. Turinici},
\it Abstract comparison principles and multivariable Gronwall-Bellman inequalities,
\rm J. Math. Anal.  Appl., 117 (1986), 100-127.



\end{thebibliography}
\end{document}